\documentstyle[amsmath,amsfonts,amssymb,amsthm,12pt]{article}
\font\sixbb=msbm6
\font\eightbb=msbm8
\font\twelvebb=msbm10 scaled 1095
\newfam\bbfam
\textfont\bbfam=\twelvebb \scriptfont\bbfam=\eightbb
                           \scriptscriptfont\bbfam=\sixbb
\def\bb{\fam\bbfam\twelvebb}
\newcommand{\Rea}{{\bb R}}

\newcommand{\Rat}{{\bb Q}}

\newtheorem{theorem}{\bf Theorem}[section]

\newtheorem{proposition}[theorem]{\bf Proposition}
\newtheorem{corollary}[theorem]{\bf Corollary}
\newcommand{\enp}{\begin{flushright} $\Box$ \end{flushright}}
\newcommand{\beq}[0]{\begin{equation}}
\newcommand{\enq}[0]{\end{equation}}
\newcommand{\th}{\tilde{{\rm H}}}

\newcommand{\cf}{{\cal F}}

\newcommand{\lk}{{\rm lk}}
\newcommand{\sd}{{\rm sd}}

\newcommand{\cs}{{\cal S}}

\newcommand{\ck}{{\cal K}}

\newcommand{\cg}{{\cal G}}

\newcommand{\dpo}{\overset{.}{D}}
\newcommand{\hell}{{\rm h}}
\newcommand{\ho}{{\rm H}}
\newcommand{\co}{{\rm C}}
\newcommand{\hv}{{\rm H}}
\newcommand{\tho}{\tilde{\rm H}}

\newcommand{\alt}{{\rm Alt}}

\newcommand{\sign}{{\rm sign}}

\newcommand{\unds}{\underline{\sigma}}
\newcommand{\tils}{\tilde{\sigma}}
\newcommand{\ler}{L}

\title{Leray Numbers of Projections and \\ a Topological Helly Type Theorem}
\author{Gil Kalai\thanks{Institute of Mathematics, Hebrew
University, Jerusalem 91904, Israel,  and Departments of Computer
Science and Mathematics, Yale University. e-mail:
kalai@math.huji.ac.il .  Research supported by ISF, BSF and NSF
grants.} \and Roy Meshulam\thanks{Department of Mathematics,
Technion, Haifa 32000, Israel. e-mail:
meshulam@math.technion.ac.il~. Research supported by the Israel
Science Foundation. }}
\begin{document}
\insert\footins{\footnotesize\rule{0pt}{\footnotesep}
\\ {\it Math Subject Classification.}  55U10,
52A35
\\ }
\maketitle
\begin{abstract}
Let $X$ be a simplicial complex on the vertex set $V$. The {\it
rational Leray number} $\ler (X)$ of $X$ is the minimal $d$ such
that $\th_i(Y;\Rat)=0$ for all induced subcomplexes $Y \subset X$
and $i \geq d$. \\
Suppose $V=\bigcup_{i=1}^m V_i$ is a partition of $V$ such that
the induced subcomplexes $X[V_i]$ are all $0$-dimensional. Let
$\pi$ denote the projection of $X$ into the $(m-1)$-simplex on the
vertex set $\{1,\ldots,m\}$ given by $\pi(v)=i$ if $v \in V_i$.
Let $r=\max\{|\pi^{-1}(\pi(x))|: x
\in |X|\}$. It is shown that
$$\ler (\pi(X)) \leq r\ler (X)+r-1~~.$$

One consequence is a topological extension of a Helly type result
of Amenta. Let $\cf$ be a family of compact sets in $\Rea^d$ such
that for any $\cf'\subset \cf$, the intersection
$\bigcap \cf'$ is either empty or contractible. \\
It is shown that if  $\cg$ is a family of sets such that for any
finite $\cg' \subset \cg$, the intersection $\bigcap \cg'$ is a
union of at most $r$ disjoint sets in $\cf$, then the Helly number
of $\cg$ is at most $r(d+1)$.
\end{abstract}

\section{Introduction}

$~~~~$ Let $\cf$ be a family of sets. The {\it Helly number} $~\hell
(\cf)$ of $\cf$ is the minimal positive integer $h$ such that if a
finite subfamily $\ck \subset \cf$ satisfies $\bigcap \ck' \neq
\emptyset$ for all $\ck' \subset \ck$ of cardinality $\leq h$, then
$\bigcap \ck \neq \emptyset$. Helly's classical theorem (1913, see
e.g. \cite{E93}) asserts that the Helly number of the family of
convex sets in $\Rea^d$ is $d+1$.

Helly's theorem and its numerous extensions are of central
importance in discrete and computational geometry (see
\cite{E93,Z97}). It is of considerable interest to understand the
role of convexity in these results, and to find suitable topological
extensions. Indeed, it is often the case that topological methods
provide a deeper understanding of the underlying combinatorics
behind Helly type theorems. Helly himself realized in 1930 (see
\cite{E93}) that in his theorem, convex sets can be replaced by
topological cells if you impose the additional requirement that all
non-empty intersections of these cells are again topological cells.
Helly's topological version of his theorem also follows from the
later nerve theorems of Borsuk, Leray and others (see below).

The following result was conjectured by Gr\"{u}nbaum  and Motzkin
\cite{MG61} , and proved by Amenta \cite{amenta96}. A family of sets
$\cg$ is an {\it $(\cf,r)$-family} if for any finite $\cg' \subset
\cg$, the intersection $\bigcap \cg'$ is a union of at most $r$
disjoint sets from $\cf$.
\begin{theorem}[Amenta]
\label{amenta} Let $\cf$ be the family of compact convex sets in
$\Rea^d$. Then for any $(\cf,r)$-family $\cg$
$$
\hell(\cg) \leq r(d+1)~~.
$$
\end{theorem}
The main motivation for the present paper was to find a topological
extension of Amenta's Theorem.

Let $X$ be a simplicial complex on the vertex set $V$. The {\it
induced} subcomplex on a subset of vertices $S \subset V$ is
$X[S]=\{ \sigma \in X: \sigma \subset S \}$. The {\it link} of a
subset $A \subset V$ is $\lk(X,A)=\{\tau \in X:\tau \cup A \in X,
\tau \cap A =\emptyset~\}~.$ The geometric realization of $X$ is
denoted by $|X|$. We identify $X$ and $|X|$ when no confusion can
arise. All homology groups considered below are with rational
coefficients, i.e. $\ho_i(X)=\ho_i(X;\Rat)$ and
$\tho_i(X)=\tho_i(X;\Rat)$.

The {\it rational Leray number} $\ler (X)$ of $X$ is the minimal $d$
such that $\th_i(Y)=0$ for all induced subcomplexes $Y \subset X$
and $i \geq d$. The Leray number can be regarded as a simple
topologically based ``complexity measure'' of $X$. Note that $\ler
(X)=0$ iff $X$ is a simplex, and $\ler (X) \leq 1$ iff $X$ is the
clique complex of a chordal graph (see \cite{W75}). It is well-known
(see e.g. \cite{KM06}) that $~\ler (X) \leq d~$ iff
$~\th_i(\lk(X,\sigma))=0~$ for all $\sigma \in X$ and $i\geq d$.
Leray numbers have also significance in commutative algebra, since
$L(X)$ is equal to the Castelnuovo-Mumford regularity of the
Stanley-Reisner ring of $X$ over $\Rat$ (see \cite{KM06}).

From now on we assume that $V_1,\ldots,V_m$ are finite disjoint
$0$-dimensional complexes, and denote their join by
$V_1*\cdots*V_m$.  Let $\Delta_{m-1}$ be the simplex on the vertex
set $[m]=\{1,\ldots,m\}$, and let $\pi$ denote the simplicial
projection from $V_1*\cdots*V_m$ onto $\Delta_{m-1}$ given by
$\pi(v)=i$ if $v \in V_i$. For a subcomplex $X \subset
V_1*\cdots*V_m$, let $r(X,\pi)=\max\{|\pi^{-1}(\pi(x))|: x \in
|X|\}$. Our main result is the following
\begin{theorem}
\label{Lproj} Let $Y=\pi(X)$ and $r=r(X,\pi)$. Then
\begin{equation}
\label{lylx} \ler (Y) \leq r\ler (X)+r-1~~.
\end{equation}
\end{theorem}
\noindent {\bf Example:} For $r \geq 1, d\geq 2$ let $m=rd$, and
consider a partition $[m]=\bigcup_{k=1}^r A_k$ with $|A_k|=d$. For
$i \in [m]$ let $V_i=\{i\}\times [r]$. Denote by $\Delta(A)$ the
simplex on vertex set $A$, with boundary $\partial \Delta(A)\simeq
S^{|A|-2}$. For $k,j \in [r]$ let $A_{kj}=A_k \times \{j\}$, and let
$$X_k= \Delta(A_{1k})* \cdots *\Delta(A_{k-1,k})
*\partial \Delta(A_{kk})*\Delta(A_{k+1,k})* \cdots
*\Delta(A_{rk})~~.$$ Let $X=\bigcup_{k=1}^r X_k$. Then $L(X)=d-1$,
and the projection $\pi:X \rightarrow \Delta_{m-1}$ satisfies
$r(X,\pi)=r$. Since $\pi(X)=\partial\Delta_{m-1}$, it follows that
$L(\pi(X))=m-1$. Hence equality is attained in (\ref{lylx}).

As mentioned earlier, Theorem \ref{Lproj} is motivated by an
application in combinatorial geometry. The {\it nerve} $N(\cf)$ of a
family of sets $\cf$, is the simplicial complex whose vertex set is
$\cf$ and whose simplices are all $\cf' \subset \cf$ such that $
\bigcap \cf' \neq \emptyset$. It is easy to see that
\begin{equation}
\label{H-L} \hell (\cf) \leq 1+ \ler (N(\cf)).
\end{equation}
A finite family $\cf$ of compact sets in some topological space is
a {\it good cover} if for any $\cf' \subset \cf$, the intersection
$\bigcap \cf'$ is either empty or contractible. If $\cf$ is a good
cover in $\Rea^d$, then  by the Nerve Lemma (see e.g. \cite{B95})
$\ler (N(\cf)) \leq d$,  hence follows the Topological Helly's
Theorem: $\hell(\cf) \leq d+1$. Theorem \ref{Lproj} implies a
similar topological generalization of Amenta's theorem.
\begin{theorem}
\label{topa} Let $\cf$ is a good cover in $\Rea^d$. Then for any
$(\cf,r)$-family $\cg$
$$\hell(\cg) \leq r(d+1)~~.
$$
\end{theorem}

The proof of Theorem \ref{Lproj} combines a vanishing theorem for
the multiple point sets of a projection, with an application of
the image computing spectral sequence due to Goryunov and Mond
\cite{GoMo93}. In Section \ref{s:icss} we describe the
Goryunov-Mond result. In Section \ref{s:hmps} we prove our main
result, Proposition \ref{hmps}, which is then used to deduce
Theorem \ref{Lproj}. The proof of Theorem \ref{topa} is given in
Section \ref{s:thelly}.

\section{The Image Computing Spectral Sequence}
\label{s:icss} $~~~$ For $X \subset V_1*\cdots*V_m$ and $k \geq 1$
define the {\it multiple point set} $M_k$ by
$$M_k=\{(x_1,\ldots,x_k) \in |X|^k~:~\pi(x_1)=\cdots=\pi(x_k)\}~~.$$

Let $W$ be a $\Rat$-vector space with an action of the symmetric
group $S_k$. Denote $\alt=\frac{1}{k!}\sum_{\sigma \in S_k}
\sign(\sigma) \sigma \in \Rat[S_k]$. Then
$$\alt\, W=\{\alt\,w: w \in W\}=$$
\begin{equation}
\label{altalt} \{w\in W:\sigma w= \sign(\sigma)w {~for~all~}
\sigma \in S_k\}~~.
\end{equation}
The natural action of $S_k$ on $M_k$ induces an action on the
rational chain complex $\co_*(M_k)$ and on the rational homology
$\ho_*(M_k)$.  The idempotence of $\alt$ implies that
\begin{equation}
\label{same} \alt\, \ho_*(M_k) \cong \ho_*(\alt\, \co (M_k))~~.
\end{equation}
The following result is due to Goryunov and Mond \cite{GoMo93} (see
also \cite{Gor95} and \cite{Hou99}).
\begin{theorem}[Goryunov and Mond]
\label{icss} Let $Y=\pi(X)$ and $r=r(X,\pi)$. Then there exists a
homology spectral sequence $\{E_{p,q}^r\}$ converging to
$\ho_*(Y)$ with
\begin{equation}
\label{epq}
E_{p,q}^1= \left\{
\begin{array}{ll}
   \alt\,\ho_q (M_{p+1}) & 0 \leq p \leq r-1, 0 \leq q \\
         0   & {\rm otherwise}
\end{array}
\right.~~
\end{equation}
\end{theorem}
\noindent {\bf Remark:} The $E^1$ terms in the original formulation
of Theorem \ref{icss} in \cite{GoMo93}, are given by
$E_{p,q}^1=\alt\,\ho_q (D^{p+1})$ where
$$D^k={\rm closure}\{(x_1,\ldots,x_k) \in |X|^k:\pi(x_1)=\cdots=\pi(x_k),
x_i \neq x_j {\rm ~for~} i \neq j \}~~.$$ The isomorphism
$$ \alt\,\ho_q (D^{p+1}) \cong \alt\,\ho_q (M_{p+1})~$$
which implies (\ref{epq}), is proved in Theorem 3.4 in \cite{Hou99}.
Indeed, as noted there, the inclusion $D^{p+1} \rightarrow M_{p+1}$
induces an isomorphism $ \alt\,\co_q (D^{p+1}) \cong \alt\,\co_q
(M_{p+1})~$ already at the alternating chains level.

\section{Homology of the Multiple Point Set}
\label{s:hmps} $~~~$

In this section we study the homology of a generalization of the
multiple point set. For subcomplexes $X_1,\ldots,X_k \subset
V_1*\cdots*V_m$, let
$$M(X_1,\ldots,X_k)= \{(x_1,\ldots,x_k) \in |X_1| \times \cdots \times |X_k|:
\pi(x_1)=\cdots =\pi(x_k)\}~.$$ In particular, if $X_1=\cdots=X_k=X$
then $M(X_1,\ldots,X_k)=M_k$.

We identify the generalized multiple point set $M(X_1,\ldots,X_k)$
with the simplicial complex whose $p$-dimensional simplices are
$\{w_{i_0},\ldots,w_{i_p}\}$, where $1 \leq i_0<\cdots<i_p \leq m$,
$w_{i_j}=(v_{i_j,1},\ldots,v_{i_j,k}) \in V_{i_j}^k$ and $\{v_{i_0
,r},\ldots,v_{i_p,r}\} \in X_r$ for all $1 \leq r \leq k$.
 The main ingredient in the proof
of Theorem \ref{Lproj} is the following
\begin{proposition}
\label{hmps} $\th_j\bigl(M(X_1,\ldots,X_k)\bigr)=0$ for $j \geq
\sum_{i=1}^k \ler (X_i)$.
\end{proposition}
The proof of Proposition \ref{hmps} depends on a spectral sequence
argument given below. We first recall some definitions. Let $K$ be
a simplicial complex. The subdivision $\sd(K)$ is the order
complex of the set of the non-empty simplices of $K$ ordered by
inclusion. For $\sigma \in K$ let $D_K(\sigma)$ denote the order
complex of the interval $[\sigma,\cdot]=\{\tau \in K ~:~ \tau
\supset \sigma \}$. Let $\dpo_K(\sigma)$ denote the order complex
of the interval $(\sigma,\cdot]=\{\tau \in K ~:~ \tau \supsetneqq
\sigma\}$. Note that $\dpo_K(\sigma)$ is isomorphic to
$\sd(\lk(K,\sigma))$ via the simplicial map $\tau \rightarrow
\tau-\sigma$. Since $D_K(\sigma)$ is contractible, it follows that
$\hv_i\bigl(D_K(\sigma),\dpo_K(\sigma)\bigr) \cong
\th_{i-1}\bigl(\lk(K,\sigma)\bigr)$
for all $i \geq 0$.  \\
For $\sigma \in V_1* \cdots* V_m$, let $\tils=\bigcup_{i\in
\pi(\sigma)} V_i$. Note that if $\sigma_2 \in X_2,\ldots,\sigma_k
\in X_k$ then there is an isomorphism
\begin{equation}
\label{observ} M(X_1,\sigma_2,\ldots,\sigma_k) \cong
X_1[\cap_{i=2}^k \tils_i]~~.
\end{equation}
For $0 \leq p \leq n=\sum_{i=2}^k \dim X_i$ let
$$
\cs'_p=\{(\sigma_2,\ldots,\sigma_k) \in X_2 \times \cdots \times
X_k~:~ \sum_{i=2}^k \dim \sigma_i \geq n-p\}~~
$$
and let $\cs_p=\cs'_p-\cs'_{p-1}$. For
$\unds=(\sigma_2,\ldots,\sigma_k) \in \cs'_p$ let
$$A_{\unds}=M(X_1,\sigma_2,\ldots,\sigma_k) \times
D_{X_2}(\sigma_2) \times \cdots \times D_{X_k}(\sigma_k)~~,$$
$$B_{\unds}=M(X_1,\sigma_2,\ldots,\sigma_k) \times \Bigl(
\bigcup_{j=2}^k D_{X_2}(\sigma_2) \times \cdots \times
\dpo_{X_j}(\sigma_j) \times \cdots \times D_{X_k}(\sigma_k)
\Bigr)~~.$$

\begin{proposition}
\label{specm} There exists a homology spectral sequence
$\{E_{p,q}^r\}$ converging to $\hv_*\bigl(M(X_1,\ldots,X_k)\bigr)$
such that
\begin{equation}
\label{spectwo}
E_{p,q}^1=\bigoplus_{\unds \in \cs_p}
\bigoplus_{\substack{i_1,\ldots,i_k \geq 0 \\
i_1+\cdots +i_k=p+q}} \hv_{i_1}\bigl(X_1[\cap_{i=2}^{k}
\tils_i]\bigr) \otimes \bigotimes_{j=2}^k
\th_{i_j-1}\bigl(\lk(X_j,\sigma_j)\bigr)
\end{equation}
 for
$~0 \leq p \leq n~$,$~0 \leq q~$, and $E_{p,q}^1=0$ otherwise.
\end{proposition}
\noindent {\bf Proof:}  For $0 \leq p \leq n$ let
$$K_p=\bigcup_{\unds \in \cs'_p} A_{\unds} \subset M(X_1,\ldots,X_k)\times
\sd(X_2) \times \cdots \times \sd(X_k).$$ Write $K=K_n$, and
consider the projection on the first coordinate $\theta:K
\rightarrow M(X_1,\ldots,X_k)$. Let $(x_1,\ldots,x_k) \in
M(X_1,\ldots,X_k)$, and let $\sigma_i$ denote the minimal simplex
in $X_i$ that contains $x_i$. Then the fiber
$$\theta^{-1}\bigl((x_1,\ldots,x_k)\bigr)= \{(x_1,\ldots,x_k)\} \times
D_{X_2}(\sigma_2) \times \cdots \times D_{X_k}(\sigma_k)$$ is a
cone, hence $K$ is homotopy equivalent to $M(X_1,\ldots,X_k)$. The
filtration $\emptyset \subset K_0 \subset \cdots \subset K_n=K$
gives rise to a homology spectral sequence $\{E_{p,q}^r\}$
converging to $\ho_*(K) \cong \ho_*(M(X_1,\ldots,X_m))$. The
$E_{p,q}^1$ terms are computed as follows. First note that
\begin{equation}
\label{preb} \bigcup_{\unds \in \cs_p} A_{\unds} \bigcap K_{p-1}
=\bigcup_{\unds \in \cs_p} B_{\unds}~~.
\end{equation}
Secondly, $\bigl(A_{\unds}-B_{\unds}\bigr) \cap A_{\unds
'}=\emptyset~$ for $~\unds \neq \unds' \in \cs_p$. Hence
\begin{equation}
\label{prea} \hv_*\Bigl(\bigcup_{\unds \in \cs_p}
A_{\unds},\bigcup_{\unds \in \cs_p} B_{\unds}\Bigr) \cong
\bigoplus_{\unds \in \cs_p} \hv_*(A_{\unds},B_{\unds})~~.
\end{equation}
 Applying excision,
(\ref{preb}),(\ref{prea}), and the K\"{u}nneth formula we obtain:
$$E_{p,q}^1=
\hv_{p+q}(K_p,K_{p-1}) \cong \hv_{p+q}\Bigl(\bigcup_{\unds \in
\cs_p} A_{\unds}, \bigcup_{\unds \in \cs_p} A_{\unds} \bigcap
K_{p-1}\Bigr) \cong
$$ $$
\hv_{p+q}\Bigl(\bigcup_{\unds \in \cs_p} A_{\unds},\bigcup_{\unds
\in \cs_p} B_{\unds}\Bigr)  \cong\bigoplus_{\unds \in \cs_p}
\hv_{p+q} (A_{\unds},B_{\unds}) \cong$$ $$ \bigoplus_{\unds \in
\cs_p}\bigoplus_{\substack{i_1,\ldots,i_k \geq 0 \\
i_1+\cdots +i_k=p+q}}
\hv_{i_1}\bigl(M(X_1,\sigma_2,\ldots,\sigma_k)\bigr) \otimes
\bigotimes_{j=2}^k
\hv_{i_j}\bigl(D_{X_j}(\sigma_j),\dpo_{X_j}(\sigma_j)\bigr)
\cong$$
$$\bigoplus_{\unds \in \cs_p}
\bigoplus_{\substack{i_1,\ldots,i_k \geq 0 \\
i_1+\cdots +i_k=p+q}} \hv_{i_1}\bigl(X_1[\cap_{i=2}^{k}
\tils_i]\bigr) \otimes \bigotimes_{j=2}^k
\th_{i_j-1}\bigl(\lk(X_j,\sigma_j)\bigr)~~.
$$
{\enp} \noindent {\bf Proof of Proposition \ref{hmps}:} If $\ler
(X_j)=0$ for all $1 \leq j \leq k$, then all the $X_j$'s are
simplicies, say $X_j=\sigma_j$. It follows that $M(X_1,\ldots,X_k)$
is isomorphic to the simplex $\bigcap_{j=1}^k \pi(\sigma_j)$ and
thus has vanishing reduced homology in all nonnegative dimensions.
Suppose then that $m=\sum_{j=1}^k \ler (X_j)>0$. Without loss of
generality we may assume that $\ler (X_1)>0$. Let $i_1,\ldots,i_k
\geq 0$ such that $\sum_{j=1}^k i_j \geq m$. Then either $i_1 \geq
\ler (X_1)$ and then $\hv_{i_1}\bigl(X_1[\cap_{i=2}^{k}
\tils_i]\bigr)=0$, or there exists a $2 \leq j \leq k$ such that
$i_j-1 \geq \ler (X_j)$ and then
$\th_{i_j-1}\bigl(\lk(X_j,\sigma_j)\bigr)=0$.  By (\ref{spectwo}) it
follows that $E_{p,q}^1=0$ if $p+q \geq m$, hence
$\th_j\bigl(M(X_1,\ldots,X_k)\bigr)=0$ for all $j \geq m$. {\enp}
\noindent {\bf Remark:} If all the $V_j$'s are singletons then
$M(X_1,\ldots,X_k)$ is isomorphic to $\bigcap_{j=1}^k X_j$. Hence
Proposition \ref{hmps} implies the following result of \cite{KM06}.
\begin{corollary}[\cite{KM06}]
\label{inter} If $X_1,\ldots,X_k$ are simplicial complexes on the
same vertex set, then
$$\ler (\bigcap_{j=1}^k X_j) \leq \sum_{j=1}^k \ler (X_j)~~.$$
{\enp}
\end{corollary}

\noindent {\bf Proof of Theorem \ref{Lproj}:} Let $Y=\pi(X)$ and
$r=r(X,\pi)$. Assuming as we may that $L(X)>0$, we have to show that
$\ho_m(Y)=0$ for $m \geq rL(X)+r-1$.
 By Theorem \ref{icss} it suffices to show that $\alt\,\ho_q
(M_{p+1})=0$ for all pairs $(p,q)$ such that $p \leq r-1$ and $p+q
\geq r\ler (X)+r-1$. Indeed,
 $p \leq r-1$ implies that $q \geq r\ler (X) \geq (p+1)\ler (X)$,
thus $\ho_q (M_{p+1})=0$ by Proposition \ref{hmps}. {\enp}

\section{A Topological Amenta Theorem}
\label{s:thelly} {\bf Proof of Theorem \ref{topa}:} Suppose
$\cg=\{G_1,\ldots,G_m\}$ is an $(\cf,r)$-family. Write
$G_i=\bigcup_{j=1}^{r_i} F_{ij}$, where $r_i \leq r$ and $F_{ij}
\bigcap F_{ij'} = \emptyset$ for $1 \leq j \neq j' \leq r_i$.  Let
$V_i=\{F_{i1},\ldots,F_{ir_i}\}$ and consider the nerve
$$X=N(\{F_{ij}:1 \leq i \leq m,~1 \leq j \leq r_i\})
\subset V_1* \cdots *V_m~~.$$ Let $\Delta_{m-1}$ be the simplex on
the vertex set $\{G_1,\ldots,G_m\}$ and let $\pi$ denote the
projection of $V_1* \cdots *V_m$ into $\Delta_{m-1}$ given by
$\pi(F_{ij})=G_i$. Then $\pi(X)=N(\cg)$. Let $y \in |N(\cg)|$ and
let $\sigma=\{G_i~:~i \in I\}$ be the minimal simplex in $N(\cg)$
such that $y \in |\sigma|$.  Then
\begin{equation}
\label{fga} |\pi^{-1}(y)|= |\{~(j_i:i \in I)~:~\bigcap_{i \in I}
F_{ij_i} \neq \emptyset~\}|~.
\end{equation}
On the other hand
\begin{equation}
\label{fgb} \bigcap_{i \in I}G_i=\bigcup_{(j_i:i \in I)}
\bigcap_{i \in I}F_{ij_i}~
\end{equation}
and the union on the right is a disjoint union. The assumption that
$\cg$ is an $(\cf,r)$ family, together with (\ref{fga}) and
(\ref{fgb}), imply that $|\pi^{-1}(y)| \leq r$ for all $y \in
|N(\cg)|$. Since $\cf$ is a good cover in $\Rea^d$, the Leray number
of the nerve satisfies $\ler(X)=\ler(N(\cf)) \leq d$. Therefore by
(\ref{H-L}) and Theorem \ref{Lproj}
$$
\hell(\cg) \leq 1+ \ler (N(\cg))=1+\ler (\pi(X)) \leq $$
$$1+ r\ler (X)+ r-1\leq r(d+1)~~.$$
{\enp}

\end{document}